\documentclass[12pt]{amsart}

\usepackage{amscd}
\usepackage{amssymb}
\usepackage{amsmath}

\theoremstyle{plain}
\newtheorem{theorem}{Theorem}[section]
\newtheorem{corollary}[theorem]{Corollary}
\newtheorem{lemma}[theorem]{Lemma}
\newtheorem{proposition}[theorem]{Proposition}

\newtheorem{definition}[theorem]{Definition}

\newtheorem{remark}[theorem]{Remark}

\begin{document}

\title [Completely positive maps and Extremal K-set]
       {Completely positive maps and
         Extremal K-set} 
\author{Hyun Ho Lee}

\address{Department Of Mathematics\\
         Purdue University\\
         West Lafayette, Indiana 47907}
\email{ylee@math.purdue.edu}
\keywords{Extremal K-set,Completely positive map,KK-theory}
\subjclass[2000]{Primary:46L0;Secondary:46L80}
\date{ May 15, 2006}
\maketitle

\begin{abstract}
   We introduce a new set $K_{e}(A,B)$ where $A$ is a commutative
   $C^{*}$-algebra, $B $ is a $C^{*}$-algebra.  It
contains $KK(A,B)$. When $A=S$ where $S$ is the
   suspension, we show there is a nice interpretation for $K_{e}(A,B)$.

\end{abstract}

\section{Introduction}
Since Kasparov invented his celebrated bivariant KK-functor in
1981 \cite{Kas81}, KK-theory was studied by several
mathematicians-Joachim Cuntz, Goerge Skandalis, Nigel Higson,
Claude Schochet, Jonathan Rosenberg during 1981-90. The power and
utility have been fully demonstrated by its applications to
geomtry, topology and recently
Elliott's Classification program since then. \\
The close connection between KK-theory and K-theory is one of main
features of KK-theory which is also the basic fact for the
Universal Coefficient Theorem (shortly UCT). This paper is a try
to pursue this point of view further along the introduction of
Brown and Pedersen 's $K_{e}$ and $E_{\infty}$\cite{BrPe}. With
the account of  the Cuntz's description of KK-group, our goal is
to find the appropriate counterpart for $K_{e}(B)$ and
$E_{\infty}(B)$ where B is a $C^{*}$-algebra. It turns out that a
slight variation of Cuntz's picture also gives us the right
candidate (See \S 4 and \S 5).\\
The plan of this papaer is as follows. After dealing with some
preliminaries on completely positive mapping in \S 2, we review
definitions of $E_{\infty}$ and $K_{e}$ and summarize the basic
facts which is necessary for our purpose in \S 3.(For our goal,
stability of these groups are essential.) In \S 4,  we define
$[\mathcal{E}(A,B)]$ and $KK_{e}(A,B)$ for $A$ which is in a small
category of commutative $C^{*}$-algebras. In contrast with the
definition KK-group, our definition shows a generality of a
variable in KK-theory comes from the restriction about morphisms.
Finally, in \S 5, we establish the connection  between
$[\mathcal{E}(A,B)]$ ,$KK_{e}(A,B)$ and $E_{\infty}(B)$,
$K_{e}(B)$ respectively which generalizes the connection
$KK^{1}(A,B)$ and $K_{1}(B)$.
\section{Completely positive map between $C^{*}$-algebras}
In this section, we show some results about completely positive
map of $C^{*}$-algebras which will be useful later.\\
For the definition and basic facts of completely positive map of
$C^{*}$-algebra, we recommend you Paulsen's Book \cite{Pl} as a
general reference.
\begin{theorem}[Stinespring's Dilation Theorem ]
Let $A$ be a (not necessarily unital) $C^{*}$-algebra and
$\mathcal{B}(\mathbb{H})$ be a space of the bounded linear
operators on Hilbert space $\mathbb{H}$. If $\phi:A \to
\mathcal{B}(\mathbb{H}) $ is completely positive map, then there
exist a Hilbert space $\mathbb{K}$, a nondegenerate
*-representation $\pi: A \to \mathcal{B}(\mathbb{K})$, and a
bounded operator $V:\mathbb{H} \to \mathbb{K}$ such that
\[
   \phi(a)=V^{*}\pi(a)V \; \mbox{for all} \; a \in A
\]

\end{theorem}
\emph{proof.}
 A general version of Stinespring Dilation theorem was suggested
 and proven by Kasparov \cite{Kas80} in Hilbert $C^{*}$-module
 setting. Most elaborate proof is found in P48-52 \cite{Lan}

\begin{lemma}\label{L:CB}
 Let $A$ and $B$ be  $C^{*}$-algebras. If $\phi:A \to B $ is
 completely positive, then $\phi $ is completely bounded and
 $\|\phi\|_{cb}=\sup_{n}\|\phi_{n}\|=\|\phi\|$. Consequently, $\phi$ is
 isometric(contractive), $\phi$ is completely
 isometric(contractive).
\end{lemma}
\emph{proof.}By the Gelfand-Naimark theorem, we can assume $\phi:
A \to \mathcal{B}(\mathbb{H})$ for some Hilbert space
$\mathbb{H}$. then by the above theorem, there exist a Hilbert
space $\mathbb{K}$, a
 *-representation $\pi: A \to \mathcal{B}(\mathbb{K})$, and a
bounded operator $V:\mathbb{H} \to \mathbb{K}$ such that
$\phi(a)=V^{*}\pi(a)V$. Since
\[
 \|\phi_{n}((a_{i,j}))\| \leq \|V^{*}\|\|(a_{i,j})\|\|V\| \;\mbox{for
 each}\; n ,
\]
$\|\phi\|_{cb}\leq \| V \|^{2}$, thus, $\phi$ is completely
bounded. But if $(e_{i})$ is an approximate unit for $A$, using
$\phi(a^{*}e_{i})\phi(e_{i}a) \leq
\|\phi(e_{i}^{2})\|\phi(a^{*}a)$, we can deduce $\|\phi\|=\sup_{i}
\| \phi(e_{i})\|$. Hence $\|\phi\|=\|V\|^{2}$. Finally, since
$e_{i}^{(n)}=e_{i}\otimes1 \in A \otimes M_{n}(\mathbb{C}) \cong
M_{n}(A)$ is an approximate unit for $M_{n}(A)$ and
$\|\phi_{n}(e_{i}^{(n)})\|=\|\phi(e_{i})\|$, using similar
inequality for $\phi_{n}$, we get
$\|\phi\|=\sup_{i}\|\phi(e_{i})\|=\sup_{i}\|\phi_{n}(e_{i}^{(n)})\|\leq
\|\phi_{n}\|$ for each $n$. So we complete the proof.

Now we are ready to prove the following proposition which will be
basic tool through this paper. We let $\tilde{A}$ be a unitization
of $A$.
\begin{proposition}
Let $A$ be non-unital, $B$ be an unital $C^{*}$-algebra and
$\phi:A \to B$ be completely positive and contractive. Then there
is a unital map $\tilde{\phi}:\tilde{A} \to B $ which is also
completely positive and extends $\phi$. Moreover, such a map is
unique.
\end{proposition}
\emph{proof.} The idea of this proof is borrowed from Huaxin Lin.
We define $\tilde{\phi}(a+\lambda)=\phi(a)+\lambda1$. Clearly,
$\tilde{\phi}$ is extends $\phi$ and unital. To show
$\tilde{\phi}$ is completely positive, we must show
$\tilde{\phi}_{n} is positive$. Note that any element $b\in
M_{n}(\tilde{A})$ can be written down as $a+s$ where $a\in
M_{n}(A)$ and $s \in M_{n}(\mathbb{C}1)$. \\
 Suppose $b$ is positive.  By considering natural quotient
 map $$M_{n}(\tilde{A}) \to \frac{M_{n}(\tilde{A})}{M_{n}(A)}$$ s
 is also positive. For any $\epsilon > 0$, if we set $s_{\epsilon}=s+\epsilon
 1$, then $s_{\epsilon}$ is invertible. Since $ a+s_{\epsilon}$ is
 positive, we have \[ -(s_{\epsilon})^{-1/2}a(s_{\epsilon})^{-1/2} \leq
 1_{M_{n}(\tilde{A})}\] Let $(e_{\lambda})$ be an approximate
 identity  for $M_{n}(A)$. Thus \[ -e_{\lambda}(s_{\epsilon})^{-1/2}a(s_{\epsilon})^{-1/2}e_{\lambda} \leq
 (e_{\lambda})^{2}.\] Since $(s_{\epsilon})^{-1/2}$ is also scalar
 matrix, $(s_{\epsilon})^{-1/2}a(s_{\epsilon})^{-1/2} \in
 M_{n}(A)$ and
 $$e_{\lambda}(s_{\epsilon})^{-1/2}a(s_{\epsilon})^{-1/2}e_{\lambda}\to
 (s_{\epsilon})^{-1/2}a(s_{\epsilon})^{-1/2}.$$ Since
 $\phi$ is completely contractive by lemma \ref{L:CB}
 \[
    \phi_{n}(-e_{\lambda}(s_{\epsilon})^{-1/2}a(s_{\epsilon})^{-1/2}e_{\lambda})\leq
    \phi_{n}(e_{\lambda})^{2})\leq 1_{M_{n}(B)}
 \]
 Observe that $\phi_{n}(sa)=s\phi_{n}(a)$ and
 $\phi_{n}(as)=\phi_{n}(a)s$ hold for any scalar matrix s. From
 this observation, we obtain
 \[
  -s_{\epsilon}^{-1/2}\phi_{n}(a)s_{\epsilon}^{-1/2}\leq 1_{M_{n}(B)}
 \] Consequently,
 \[
   \phi_{n}(a)+ S_{\epsilon}\geq 0
 \]  Thus $$\tilde{\phi}(a+s)=\phi_{n}(a)+s \geq -\epsilon
 1_{M_{n}(B)}$$ Now let $\epsilon \to 0$. we get $$\tilde{\phi}(a+s)\geq 0.
 $$
 So we have shown $\tilde{\phi}_{n}$ is positive.Uniqueness part is obvious.

\begin{lemma}\label{L:S}
   Let $A$ be a commutative $C^{*}$-algebra and $B$ a
   $C^{*}$-algebra. If $\phi:A \mapsto B$ is positive, then $\phi$
   is completely positive.
\end{lemma}
\emph{proof.} (Unital case) If $A$ is a unital commutative
$C^{*}$algebra, then we can assume that $A=C(X)$ where $X$ is a
compact Hausdorff space. Let $\epsilon > 0$ be given and $P(x)$ be
positive in $\mathbb{M}_{n}(C(X))$. We must prove $\phi_{n}(P)$ is
positive. Using a partition of  unity $\{u_{l}(x)\}$ subordinate
to a covering $\{O_{l}\}$ such that $\|P(x) - P(x_{l})\| <
\epsilon$ for $x\in O_{l}$ and positive matrices
$P(x_{l})=P_{l}=(P_{i,j})$, we have
\[
   \|P(x)-\sum_{l}u_{l}(x)P_{l}\| < \epsilon
\]
But $\phi_{n}(u_{l}P_{l})=(\phi(u_{l})p_{i,j})$ which is postive
in $\mathbb{M}_{n}(C(X))$. Since $M_{n}(B)^{+}$ is closed set,
$\phi_{n}(P)$ is positive.\\
(Non-unital case) If A is non-unital, we extend $\phi$ to
$\tilde{A}$ the unitization of A. Define $\tilde{\phi}:\tilde{A}
\to \tilde{B}$ by $\tilde{\phi}(a+\lambda1)=\phi(a)+\lambda
\|\phi\|1$. Note that if $a+\lambda1$ is positive, then $\lambda
\geq 0$. From this, $\tilde{\phi}$ is also positive. Hence
$\tilde{\phi}$ is completely positive as we have seen above.
Consequently, Restriction of $\tilde{\phi}$ to A which is $\phi$
is also completely positive.

\begin{lemma}\label{L:Observation}
Let $B$ be a unital $C^{*}$-algebra and $A$ be a $C^{*}$-
algebra(not necessarily unital) $\phi:A \to B$ be a completely
positive contractive map. then
 \begin{enumerate}
\item $ \{ a \in A |\,\phi(a)\phi(a^{*})=\phi(aa^{*}) \} = \{a\in
A | \,
 \phi(a)\phi(b)=\phi(ab) \,\mbox{for all}\,\, b \in A  \}$ is a
subalgebra of A and $\phi$ is a homomorphism when it is restricted
to this set. \item  $ \{ a \in A|\,\phi(a^{*})\phi(a)=\phi(a^{*}a)
\} = \{a\in A | \, \phi(b)\phi(a)=\phi(ba) \,\mbox{for all}\,\, b
\in A  \}$ is a subalgebra of A and $\phi$ is a homomorphism when
it is restricted to this set.
 \item $ \{ a \in A |\,
\phi(a)\phi(a^{*})=\phi(aa^{*}) \&
\phi(a^{*}a)=\phi(a^{*})\phi(a)\} $ \\ $ =\{a\in A | \,
\phi(a)\phi(b)=\phi(ab) \& \phi(ba)=\phi(b)\phi(a)\, \mbox{for
all} \,\, b \in A \}$ is a subalgebra of A and $\phi$ is a
*-homomorphism when it is restricted to this set.
 \end{enumerate}
\end{lemma}

\emph{proof.}
we prove (1). The proofs of (2) and (3) are similar.\\
We may assume that $\phi:A \mapsto \mathcal{B}(\mathbb{H})$ where
$\mathbb{H}$ is a Hilbert space. By the theorem, there is a
Hilbert space $\mathbb{K}$ containing $\mathbb{H}$, $V \in
\mathcal{B}(\mathbb{K})$ with  $\| V \| \leq 1$ and $\pi:A \mapsto
\mathcal{B}(\mathbb{K})$ is a *-representation of $A$ such that
$\phi(a)=V^{*}\pi(a)V$ for all $a$.\\
Let $a$ belong to the set on the left. then since
$\phi(a)\phi(a^{*})=\phi(aa^{*})$ holds, we have
$V^{*}\pi(a)(1-VV^{*})\pi(a)V=0$. Note that $\| V \| \leq 1$. This
implies that $1-VV^{*}$ is positive. Hence
\[
V^{*}\pi(a)(1-VV^{*})^{1/2}(1-VV^{*})^{1/2}\pi(a^{*})V=0
\]Consequently, $V^{*}\pi(a)(1-VV^{*})^{1/2}=0$.Then
\begin{eqnarray*}
 V^{*}\pi(a)(1-VV^{*})^{1/2}(1-VV^{*})^{1/2}\pi(b)V &=& 0 \,\mbox{for all}\,\, b \in A \\
     V^{*}\pi(a)(1-VV^{*})\pi(b)V &=& 0 \,\mbox{for all}\,\, b \in A \\
     \therefore V^{*}\pi(a)VV^{*}\pi(b)V &=& V^{*}\pi(ab)V
\end{eqnarray*}
So we have shown that $a$ is the element of the set on the right.

\begin{remark}
  (i) We call the set in (1) \emph{left multiplicative domain} for $\phi$, the
 set in (2) \emph{right multiplicative domain} for $\phi$ and the set in
 (3) \emph{multiplicative domain} for $\phi$.\\
  (ii) There is a more general virsion of above lemma. See the theorem 3.18 in \cite{Pl}.
\end{remark}

\begin{corollary}\label{Observation2}
Let $\phi: S \to \mathcal{M}(B\otimes \mathcal{K})$ be completely
positive contractive map. If $\phi(f)+1$ is a unitary, then $\phi$
is *-homomorphism.
\end{corollary}

\emph{proof}
 If $\phi(f)+1$ is a unitary, then $\phi(f)\phi(f^{*})=\phi(ff^{*}) \&
 \phi(f^{*})\phi(f)=\phi(f^{*}f)$ hold. Note that $f$ is a generator of
 $S$. Hence $ \{ a \in S |\,
\phi(a)\phi(a^{*})=\phi(aa^{*}) \&
\phi(a^{*}a)=\phi(a^{*})\phi(a)\} $ is nonempty and is S itself.
$\therefore  \phi$ is a *-homomorphism by
lemma\ref{L:Observation}.

\begin{proposition}\label{P:CP}
  Let $A$ be a $C^{*}$-algebra(not necessarily unital), $B$ be a
  unital $C^{*}$-algebra. Let $\phi:A \to B $ be a positive map.
  Assume $I_{+},I_{-}$ are centrally orthogonal ideals in $B$. Let $\pi^{+},\pi^{-}$ be
  the natural quotient maps from $B$ onto $B/I_{+},B/I_{-}$ respectively.
  If
  \begin{eqnarray*}
    \phi^{+}= \pi^{+}\circ \phi: A \to B \to B/I_{+} \\
    \phi^{-}=\pi^{-} \circ \phi: A \to B \to B/I_{-}
  \end{eqnarray*}
  are completely positive contractive, then $\phi:A \to B$ is
  completely positive and contractive.
\end{proposition}
\emph{proof.} Consider the *-homomorphism $\pi^{+}\oplus\pi^{-}:B
\to B/I_{+}\oplus B/I_{-}$. Then this map is injective since
$I^{+} \cap I^{-}=0$. Since the map $\phi^{+}\oplus\phi^{-}:A \to
B/I_{+}\oplus B/I_{-}$ is completely positive by the assumtion,
$\phi$ is completely positive. In fact, $B \cong \{ x\oplus y \in
B/I_{+}\oplus B/I_{-} \mid \,\mbox{Image of x} = \mbox{Image of y}
\; \mbox{in} \, B/ I^{+}+I^{-} \}$. Hence $\|\phi(a) \|=\max\{\|
\phi^{+}(a)\| ,\| \phi^{-}(a) \| \} \leq \|a\|$. So we have shown
$\phi$ is also contractive.

\section{Stable extremal class and $K_{e}(-)$}
 In this section, we summarize the definitions of $E_{\infty}(-)$
 and $K_{e}(-)$ and basic results from \cite{BrPe}.
 Throughout this section $A$ will denote a unital $C^{*}$-algebra
 and $\mathcal{E}(A)$ the set of extreme points in the unit ball
 $A_{1}$ of $A$,that is, the patial isometries $v$ such that  $(1-v^{*}v)A(1-vv^{*})=
 0$. The centrally orthogonal projections $p_{+} = 1-v^{*}v$ and $p_{-} =
 1-vv^{*}$ will be referred to as the \emph{defect projections} of $v$, and the two orthogonal
  closed ideals $I_{+}$ and $I_{-}$ generated by these
 projections will be known as the \emph{defect ideals} of $v$.

\begin{theorem}
Given $\epsilon >0$ there is a $\delta > 0$, such that for any
pair $v, w$ in $\mathcal{E}(A)$ with $\|v-w \|\leq \delta$, there
are unitaries $u_{1}$ and $U_{2}$ in $A$ with $v=u_{1}wu_{2}$ and
$\|1-u_{i}\|\leq \epsilon$ for $i=1,2$.
\end{theorem}
\emph{proof}
 This is the theorem 2.1 in \cite{BrPe}.

\begin{corollary}
Two elements $v$ and $w$ in $\mathcal{E}(A)$ are homotopic if and
only if $w=u_{1}wu_{2}$ for some unitaries $u_{1}, u_{2}$ in
$\mathcal{U}_{0}(A)$, the connected component of the unitary group
$\mathcal{U}(A)$ containing $\bold{1}$.
\end{corollary}
\emph{proof.} See the corollary 2.3. in \cite{BrPe}.

\begin{proposition}\label{P:Composability}
Let $v$ and $w$ be extremal partial isometries in $A$, and consider
the defect projections $p_{+}=1-v^{*}v, p_{-}=1-vv^{*}$, and
$q_{+}=1-w^{*}w,q_{-}=1-ww^{*}$; and the corresponding defect
ideals $I_{+},I_{-}$ and $J_{+},J_{-}$. Then the following are
equivalent.
\begin{enumerate}
 \item $p_{+}Aq_{-} = p_{-}Aq_{+}=\{0 \}$,
 \item $I_{+} \cap J_{-}= I_{-}\cap J_{+}= \{0 \}$,
 \item $vw\in \mathcal{E}(A)$ and $ wv \in \mathcal{E}(A)$,
 \item $ \left(\begin{matrix}
            v & 0 \\
            0 & w
       \end{matrix}\right) \in \mathcal{E}(M_{2}(A))$
\end{enumerate}
\end{proposition}
\emph{proof.}
 This is the proposition 2.5 in \cite{BrPe}.

When the conditions in proposition\ref{P:Composability} are
satisfied, we say that $v$ and $w$ are \emph{composable}. We see
the unitary elements are composable to all other elements in
$\mathcal{E}(A)$, and also that they are the only elements
composable their adjoints.
\begin{corollary}\label{C:composabilityhomotopy}
Suppose $v$ and $w$ are composable. If $v_{1}$ is homotopic to $v$
and $w_{1}$ is homotopic to $w$, then $v_{1}$ and $w_{1}$ is also
composable.
\end{corollary}
\emph{proof}
 It is easy to check that defect ideals of $v_{1}$ and $w_{1}$ are
 exactly the defect ideals of $v$ and $w$ respectively. From this,
 the conclusion is straightforward.

\begin{proposition}
When $v$ and $w$ are composable, two defect ideals for vw are
precisely $I_{+}+I_{-}$ and $I_{-}+J_{-}$.
\end{proposition}
\emph{proof}
 Note that $P_{+}=1-(vw)^{*}vw=q_{+} +w^{*}p_{+}w $. Defect ideal which
 is generated by $P_{+}=1-(vw)^{*}vw$ is contained in $I_{+}+I_{-}$.
 Conversely,from $wP_{+}w^{*}=p_{+} and P_{+}q_{+}=q_{+}$,
 $I_{+}+I_{-}$is also contained the ideal generated by $P_{+}$.
 The other case is proved similarly.

\begin{definition}
Let $A$ be a unital $C^{*}$-algebra and for each n, consider the
set $[\mathcal{E}(\mathbb{M}_{n}(A))]$ of homotopy classes of
extreme partial isometires in the algebra $\mathbb{M}_{n}(A)$. The
embeddings  \[ \iota_{mn}: v \hookrightarrow \left( \begin{matrix} v & 0 \\
                    0 & 1_{m-n}  \end{matrix} \right) \]
  of $\mathcal{E}(\mathbb{M}_{n}(A))$ into
  $\mathcal{E}(\mathbb{M}_{m}(A))$, for $1\leq n < m$, evidently
  respect homotopy, we define
  $[\mathcal{E}_{\infty}(A)]=\varinjlim[\mathcal{E}(\mathbb{M}_{n}(A))]$ in
  complete analogy with the definition of $K_{1}(A)$. We shall
  refer this set as the set of \emph{stable extremal classes} for
  $A$ and we shall denote by $[v]$ its homotopy class in
  $[\mathcal{E}_{\infty}(A)]$.
\end{definition}

\begin{proposition}\label{P:ESC}
Let $A$ be a $C^{*}$-algebra. Then
$$ [\mathcal{E}_{\infty}(A)]\cong
[\mathcal{E}((A\otimes \mathcal{K} )^{\sim})]$$
\end{proposition}
\emph{proof.}
 Note that $\iota_{mn}$ induce the isomorphism betweem
 $[\mathcal{E}_{\infty}(\mathbb{M}_{n}(A))]$ and
 $[\mathcal{E}_{\infty}(\mathbb{M}_{m}(A))]$ for $ 1\leq n < m $.
 In fact, any 'corner' embeddings $id\otimes e_{ii}$ are all
 homotopic. \\
 If $A=\varinjlim A_{i}$ is an s.e.p.p. inductive limit of
 $C^{*}$-algebras(i.e. every connecting morphism $\iota_{i,j}$ is stably extremal preserving
 maps), then
 $[\mathcal{E}_{\infty}(A)]=\varinjlim[\mathcal{E}_{\infty}(A_{i})]$
 if each $\iota_{i,j}$ is injective.( See P218 in \cite{BrPe}).\\
 Then by combining these two observations, we see that
 $[\mathcal{E}_{\infty}(A\otimes\mathcal{K})]=[\mathcal{E}_{\infty}(A)]
 $.\\ Since $[\mathcal{E}(M_{n}(A\otimes\mathcal{K}))]
 \cong[\mathcal{E}(A\otimes\mathcal{K}\otimes\mathbb{M}_{n})]\cong[\mathcal{E}(A\otimes\mathcal{K})]$,
$[\mathcal{E}_{\infty}(A\otimes\mathcal{K})]\cong
[\mathcal{E}(A\otimes\mathcal{K})]=[\mathcal{E}((A\otimes\mathcal{K})^{\sim})]$

\begin{proposition}
Let $A$ be a $C^{*}$-algebra. If $id\otimes e_{11}: A \mapsto
A\otimes \mathcal{K}$ is the corner embedding, then
$K_{e}(id\otimes e_{11})$ is an isomorphism between $K_{e}(A)$ and
$K_{e}(A\otimes \mathcal{K})$.
\end{proposition}
\emph{proof.} The argument is almost identical to proposition
\ref{P:ESC}. See page 219 4.12 (iv) of \cite{BrPe}.

\section{The set  $[\mathcal{E}_{\infty}(-,-)]$ and $KK_{e}(-,-)$ }
Based upon Brown and Pedersen's work \cite{BrPe}, we define more
general functor $[\mathcal{E}_{\infty}(A,B)]$ where $A$ is a
$C^{*}$-algebra generated by $u-1$ where $u$ is an unitary on a
Hilbert space and $B$ is $\sigma$-unital $C^{*}$-algebra.
\begin{definition}\label{D:cycle}
Let $\mathcal{E}(A,B)$ be the set of the pairs
$(\phi_{+},\phi_{-})$\,s.t.
\begin{enumerate}
\item $\phi_{+}:A \to \mathcal{M}(B\otimes\mathcal{K}) $
*-homomorphism,
      $ \phi_{-}: A \to \mathcal{M}(B\otimes\mathcal{K})$ completely positive contractive s.t.
    \begin{enumerate}
      \item $\tilde{\phi_{-}}(u)$ is extremal partial isometry for the unitary element $u \in \tilde{A}$ \\
      \item $ \tilde{\phi_{-}}(u^{n})=(\tilde{\phi_{-}}(u))^{n}$
            $\tilde{\phi_{-}}((u^{*})^{n})=(\tilde{\phi_{-}}(u^{*}))^{n}$ for the unitary element $u \in \tilde{A}$
     \end{enumerate}
\item $\tilde{\phi_{+}}(u) - \tilde{\phi_{-}}(u) \in
B\otimes\mathcal{K}$ thus $(\tilde{\phi_{+}}(u))^{*}
\tilde{\phi_{-}}(u) \in 1+B\otimes\mathcal{K} $ \\
\end{enumerate}
\end{definition}

We call such a pair $(\phi_{+},\phi_{-})$ as \emph{a generalized
extremal cycle} and  two generalized extremal cycles
$(\phi_{+},\phi_{-}),(\psi_{+},\psi_{-})$ are \emph{homotopic}
   when there is a generalized extremal cycle from $A$ to
$\mathcal{M}([0,1]B\otimes \mathcal{K})$ i.e.
$(\lambda_{+}^{t},\lambda_{-}^{t})\in
   \mathcal{E}(A,B), t \in [0,1] $ s.t.
   \begin{enumerate}
   \item the maps $t \to \tilde{\lambda_{+}^{t}}(u)$ and $ t \to
\tilde{\lambda_{-}^{t}}(u)$ from $[0,1]$ to
$\mathcal{M}(B\otimes\mathcal{K})$
         are strictly continuous.
   \item $\tilde{\lambda_{+}^{t}}(u)- \tilde{\lambda_{-}^{t}}(u) \in B\otimes\mathcal{K}$ for each
$t$. and the map $ t \to \tilde{\lambda_{+}^{t}}(u)-
\tilde{\lambda_{-}^{t}}(u)$ from $[0,1]$ to $B\otimes\mathcal{K}$
is norm continuous and thus the map $ t \to
(\tilde{\lambda_{+}^{t}}(u))^{*}\tilde{\lambda_{-}^{t}}(u)$ from
$[0,1]$ to $1+B\otimes\mathcal{K}$ is norm continuous
   \end{enumerate}
   We write $(\phi_{+},\phi_{-})\sim (\psi_{+},\psi_{-})$ in this
   case.
  \begin{definition}
   We let $[\mathcal{E}_{\infty}(A,B)]\overset{\text{def}}{=}\mathcal{E}(A,B)/\sim$ denote the
homotopy classes of generalized extremal cycles. The homotopy
classes in $E(A,B)$ represented by $(\phi_{+},\phi_{-})\in
\mathcal{E}(A,B)$ is denoted by $[\phi_{+},\phi_{-}]$.
  \end{definition}
We shall refer to $[\mathcal{E}_{\infty}(A,B)]$ as the set of
\emph{stable extremal classes} of $A$ and $B$.
  \begin{proposition}\label{P:Subset}
   $KK(A,B) \subseteq [\mathcal{E}_{\infty}(A,B)] $.
  \end{proposition}
  \emph{proof.}
   Since $KK(A,B)=\{[\phi_{+},\phi_{-}] \mid \phi_{+}(a)-\phi_{-}(a)\in
   B\otimes\mathcal{K}\;\mbox{for each}\; a\in A \}$,
   where $\phi_{+},\phi_{-}\in \mbox{Hom}(A,\mathcal{M}(B\otimes\mathcal{K}))$(See Chapter4 in \cite{JK}), their unital extensions
   $\tilde{\phi_{+}},\tilde{\phi_{-}}$ are *-homomorphisms. In
   particular, $\tilde{\phi_{-}}$ satisfies the conditions  (a)and (b) in definition\ref{D:cycle}.
   Since unitaries are composable to any extremal partial
   isometry, $\tilde{\phi_{+}}(u)^{*}\tilde{\phi_{-}}(u),
\tilde{\phi_{+}}(u)\tilde{\phi_{-}}(u)^{*}$ is of
   the form $1+B\otimes\mathcal{K}$ because
   $\tilde{\phi_{+}}(u)-\tilde{\phi_{-}}(u) \in
   B\otimes\mathcal{K}$. i.e. $(\phi_{+},\phi_{-})$ is in
   $\mathcal{E}(A,B)$. We can apply the same argument to a homotopy
   between two KK-cycles so that it is well-defined.

We can define a partial addition for manageable extremal classes.
 \begin{definition}
 We say two generalized extremal cycles $(\phi_{+},\phi_{-})$ and $(\psi_{+},\psi_{-})$ are
 \emph{composable} if $\tilde{\phi_{-}}(u)$ and $\tilde{\psi_{-}}(u)$ are
composable in $\mathcal{E}(\mathcal{M}(B\otimes\mathcal{K}))$
  for $u$.
\end{definition}

  Suppose two generalized extremal cycles $(\phi_{+},\phi_{-})$ and $(\psi_{+},\psi_{-})$ are
 \emph{composable}. Then we can check the following facts.
 \begin{enumerate}
  \item $ \Theta_{B}\circ \left [
                                \begin{matrix}
                                \tilde{\phi_{+}} & 0     \\
                                        0        & \tilde{\psi_{+}}
                               \end{matrix}
                           \right ]$ is *-homomorphism ,
         $ \Theta_{B}\circ \left[
                                \begin{matrix}
                                \tilde{\phi_{-}} & 0     \\
                                        0        & \tilde{\psi_{-}}
                               \end{matrix}
                           \right]$ is also completely positive
                           map which satisfies the conditions
                           $(a),(b)$ of  definition \ref{D:cycle}.
   \item     $\left(\Theta_{B}\circ \left[
                                \begin{matrix}
                                \tilde{\phi_{+}} & 0     \\
                                        0        & \tilde{\psi_{+}}
                               \end{matrix}
                           \right]\right)(u) -  \left(\Theta_{B}\circ \left[
                                \begin{matrix}
                                \tilde{\phi_{-}} & 0     \\
                                        0        & \tilde{\psi_{-}}
                               \end{matrix}
                           \right]\right)(u) \in B\otimes \mathcal{K}  \\
                              \left( \Theta_{B}\circ \left[
                                \begin{matrix}
                                \tilde{\phi_{+}} & 0     \\
                                        0        &
                                        \tilde{\psi_{+}}
                               \end{matrix}
                           \right] \right )^{*}(u) \, \left( \Theta_{B}\circ \left[
                                \begin{matrix}
                                \tilde{\phi_{-}}(u) & 0     \\
                                        0        &
                                        \tilde{\psi_{-}}(u)
                               \end{matrix}
                           \right] \right)(u) \\
                            $ are of the form $1+
                           B\otimes\mathcal{K}$ \; \mbox{for}\;
                           $u$.
\end{enumerate}
where $\Theta_{B}:M_{2}(\mathcal{M}(B\otimes\mathcal{K})) \to
\mathcal{M}(B\otimes\mathcal{K})$ is an inner *-isomorphism.\\
Now we can define an addition between two \emph{composable}
generalized extremal classes by
 \[
   [\phi_{+},\phi_{-}] + [\psi_{+},\psi_{-}]= \left[ \Theta_{B}\circ\left[
                                \begin{matrix}
                                \phi_{+} & 0     \\
                                        0        & \psi_{+}
                               \end{matrix}
                           \right], \Theta_{B}\circ \left[
                                \begin{matrix}
                                \phi_{-} & 0     \\
                                        0        & \psi_{-}
                               \end{matrix}
                           \right]
                           \right]
\]

It is easy to check each element in $ KK(A,B)$ is
\emph{composable} to an element in $E(A,B)$. This implies our
definition of \emph{composability} follows the same spirit of the
definition of \emph{composability} in \cite{BrPe}. Next result
consolidates our definition of \emph{composability} is exactly
analogous to and extended notion of the definition of
\emph{composability} of extremal partial isometries.
\begin{proposition}
The addition is associative whenever possible.
\end{proposition}
\emph{proof} If $\alpha, \beta, $ and $\gamma$ are elements in
$E(A,B)$ such that $\alpha$ is composable with $\beta$ and $\alpha
+ \beta$ is composable with $\gamma$, then $\beta$ is composable
with $\gamma$ and $\alpha$ is composable with $\beta + \gamma$.
This follows by observing that defect ideal for a sum of elements
is the sum of the defect ideals for the summands. Then using
rotational homotopies finally we have $(\alpha+\beta)+\gamma =
\alpha + (\beta + \gamma) $(See Lemma 1.3.12 in \cite{JK}).

We summarize our observations in the following theorem.
\begin{theorem}
For a non-unital $C^{*}$-algebra $A$ which is generated by $u-1$
where $u$ is an unitary on a Hilbert space and $\sigma$-unital
$C^{*}$-algebra $B$  the set of extremal classes of $A$ and $B$
$[\mathcal{E}_{\infty}(A,B)]$ is a set with a partially defined
addition between composable elements. There is a natural embedding
$KK(A,B)\subset [\mathcal{E}_{\infty}(A,B)]$, and addtion in
$[\mathcal{E}_{\infty}(A,B)]$ extends the addition in $KK(A,B)$.
\end{theorem}
 As Brown and Pedersen have defined a coarser equivalence
relation than homotopy\cite{BrPe} , we introduce a coarser
equivalence relation than homotopy to $\mathcal{E}(A,B)$.
\begin{definition}
For any two elements $\alpha=[\phi_{+},\phi_{-}]$ and
$\beta=[\psi_{+}, \psi_{-}]$ in $E(A,B)$, we define $\alpha
\approx \beta $ in ${E}(A,B)$ if there is an element $(\tau_{+},
\tau_{-})\in \mathcal{E}(A,B)$ s.t.
\begin{enumerate}
 \item $\tilde{\tau_{+}}(u)$ has smaller defects than
$\tilde{\phi_{+}}(u),\tilde{\psi_{+}}(u)$.
  \item $\tilde{\tau_{-}}(u)$ has smaller defects than
$\tilde{\phi_{-}}(u),\tilde{\psi_{-}}(u)$.
 \item $[\phi_{+},\phi_{-}]+[\tau_{+},\tau_{-}]= [\psi_{+},\psi_{-}]+[\tau_{+},\tau_{-}]
 $ or, \\ \[\left(  \Theta_{B}\circ \left[
                                \begin{matrix}
                                \phi_{+} & 0     \\
                                        0        & \tau_{+}
                               \end{matrix}
                           \right], \Theta_{B}\circ \left[
                                \begin{matrix}
                                \phi_{-} & 0     \\
                                        0        & \tau_{-}
                               \end{matrix}
                           \right] \right) \sim \left( \Theta_{B}\circ\left[
                                \begin{matrix}
                                \psi_{+} & 0     \\
                                        0        & \tau_{+}
                               \end{matrix}
                           \right], \Theta_{B}\circ \left[
                                \begin{matrix}
                                \psi_{-} & 0     \\
                                        0        & \tau_{-}
                               \end{matrix}
                           \right] \right)\]
\end{enumerate}
Evidently, we may assume that $\alpha$
and $\beta$ had the same defect ideals. i.e. $\tilde{\phi_{+}}(u),
\tilde{\psi_{+}}(u)$ had the same defect ideals  and $
\tilde{\phi_{-}}(u), \tilde{\psi_{-}}(u)$ had the same ideals for
each $u \in \mathcal{E}(A)$.\\
\end{definition}
 To verify this is an equivalence relation, we only prove
  transitivity part.\\
  Now we let $\alpha=[\phi^{1}_{+},\phi^{1}_{-}] \approx \beta=[\phi^{2}_{+},\phi^{2}_{-}]$ and
  $\beta= [\phi^{2}_{+},\phi^{2}_{-}]\approx
  \gamma=[\phi^{3}_{+},\phi^{3}_{-}]$. There is
  $(\mu_{+},\mu_{-})$ and $(\nu_{+},\nu_{-})$ in $\mathcal{E}(A)$ s.t.
\begin{subequations}
  \begin{equation}\label{E:1}
    \left(  \Theta_{B}\circ \left[
                                \begin{matrix}
                                \phi^{1}_{+} & 0     \\
                                        0        & \mu_{+}
                               \end{matrix}
                           \right], \Theta_{B}\circ \left[
                                \begin{matrix}
                                \phi^{1}_{-} & 0     \\
                                        0        & \mu_{-}
                               \end{matrix}
                           \right] \right) \sim \left( \Theta_{B}\circ\left[
                                \begin{matrix}
                                \phi^{2}_{+} & 0     \\
                                        0        & \mu_{+}
                               \end{matrix}
                           \right], \Theta_{B}\circ \left[
                                \begin{matrix}
                                \phi^{2}_{-} & 0     \\
                                        0        & \mu_{-}
                               \end{matrix}
                           \right] \right)
  \end{equation} and
  \begin{equation}\label{E:2}
    \left(  \Theta_{B}\circ \left[
                                \begin{matrix}
                                \phi^{2}_{+} & 0     \\
                                        0        & \nu_{+}
                               \end{matrix}
                           \right], \Theta_{B}\circ \left[
                                \begin{matrix}
                                \phi^{2}_{-} & 0     \\
                                        0        & \nu_{-}
                               \end{matrix}
                           \right] \right) \sim \left( \Theta_{B}\circ\left[
                                \begin{matrix}
                                \phi^{3}_{+} & 0     \\
                                        0        & \nu_{+}
                               \end{matrix}
                           \right], \Theta_{B}\circ \left[
                                \begin{matrix}
                                \phi^{3}_{-} & 0     \\
                                        0        & \nu_{-}
                               \end{matrix}
                           \right] \right)
  \end{equation}
\end{subequations}
Note that $(\mu_{+},\mu_{-})$ and $(\nu_{+},\nu_{-})$ are
composable , $ \Theta_{B}\circ\left[
                                \begin{matrix}
                                \mu_{+}(u) & 0     \\
                                        0        & \nu_{+}(u)
                               \end{matrix}
                           \right]$ has smaller defects than
                           $\phi^{1}_{+}(u)$ and $ \Theta_{B}\circ\left[
                                \begin{matrix}
                                \mu_{-}(u) & 0     \\
                                        0        & \nu_{-}(u)
                               \end{matrix}
                           \right]$ has smaller defects than
                           $\phi^{1}_{-}(u)$. \\

Then
\begin{align*}
 \alpha+[\mu_{+},\mu_{-}]+[\nu_{+},\nu_{-}]&= \beta+[\mu_{+},\mu_{-}]+[\nu_{+},\nu_{-}]\; \mbox{by}\; (\ref{E:1})\\
                                   &=\beta+[\nu_{+},\nu_{-}]+[\mu_{+},\mu_{-}]\\
                                   &=\gamma+[\nu_{+},\nu_{-}]+[\mu_{+},\mu_{-}]\; \mbox{by} \; (\ref{E:2})\\
                                   &=\gamma+[\mu_{+},\mu_{-}]+[\nu_{+},\nu_{-}]
\end{align*}

We define
\[
     KK_{e}(A,B)=[\mathcal{E}_{\infty}(A,B)]/\approx
\]
and we shall refer to $ KK_{e}(A,B)$ as the \emph{extremal KK-set}
of A and B. \\
If $\alpha$ and $\beta$ are elements in $KK(A,B)$ and $\alpha
\approx \beta$,  then only choice of $[\tau_{+},\tau_{-}]$ with
smaller defect ideals  is another element in $KK(A,B)$, whence
$\alpha=\beta$. We therefore have a natural embedding of $KK(A,B)$
into $KK_{e}(A,B)$. The natural class map $\kappa_{e}$ from
$[\mathcal{E}_{\infty}(A,B)]$ onto $KK_{e}(A,B)$ respects
composability and addition.( Here composable classes in
$KK_{e}(A,B)$ meaning that one, hence any, pair of representatives
in $[\mathcal{E}_{\infty}(A,B)]$ are composable.)\\
We summarize our observations in the following theorem.
\begin{theorem}
For a non-unital $C^{*}$-algebra $A$ which is generated by $u-1$
where $u$ is an unitary on a Hilbert space and $\sigma$-unital
$C^{*}$-algebra $B$ the extremal KK-set of $A$ and $B$ is the set
with a partially defined addition between composable elements.
There is a natural embedding $KK(A,B) \subset KK_{e}(A,B)$, and
the addition in $KK_{e}(A,B)$ extends the addition $KK(A,B)$.
There is a natural map $\kappa_{e}:[\mathcal{E}_{\infty}(A,B)] \to
KK_{e}(A,B)$ which is surjective and restricts to an isomorphism
on $KK(A,B)$.
\end{theorem}
\section{Application: special case $[\mathcal{E}_{\infty}(S,-)]$  }
We begin this section by observing the following proposition. Let
$z(t)=e^{2 \pi t i} \in C(\mathbb{T})$, $f(t)=e^{2\pi t i}-1 \in S
$.
\begin{proposition}\label{P:Suspension}
  Let $\phi:S \to \mathcal{M}(B\otimes \mathcal{K})$ be a
  completely positive contractive map and $\tilde{\phi}:C(\mathbb{T}) \to
\mathcal{M}(B\otimes \mathcal{K})  $ is the
  unital extension of
  $\phi$. Then the following are
  equivalent.\\
\begin{enumerate}
  \item\quad $\phi(f \bar{f})-\phi(f)\phi(\bar{f})\mathcal{M}(B\otimes
\mathcal{K})\phi(\bar{f} f)-\phi(\bar{f})\phi(f)=0
  $
  \item\quad $\phi(f)+1 $ is a extremal partial isometry.
  \item \quad $\tilde{\phi}(z)$ is a extremal partial isometry.

\end{enumerate}
\end{proposition}

\emph{proof.} $ (1)\Rightarrow(2)\Leftrightarrow(3)\Leftarrow(4) $
is
easy.\\
  For $(1)\Leftarrow(2)$,\; let $\phi(f)+1 = V$. Note that
  $ f\bar{f}+f+\bar{f}=0 $. Then it is easy to check
  $1-VV^{*}=\phi(f\bar{f})-\phi(f)\phi(\bar{f})$.
  Similarly, $1-V^{*}V=\phi(\bar{f}f)-\phi(\bar{f})\phi(f)$.\\

 \begin{proposition}
 Given a completely positive contractive map $\phi: C(\mathbb{T}) \to B$ a $C^{*}$- algebra,
 if ${\phi}(z)$ is an extremal partial isometry and
  $\phi(e^{ih})=e^{i\phi(h)}$ for any self-adjoint $h$ in $C(\mathbb{T})$,
  then $\tilde{\phi}$ is an extremal preserving map.
  \end{proposition}
 \emph{proof.}  we assume $\tilde{\phi}(z)$ is an
  extremal partial isometry, say $v$.\\
   Consider $\phi^{+}= \pi^{+}\circ \phi: C(\mathbb{T}) \to \
  \mathcal{M}(B\otimes\mathcal{K}) \to
  \mathcal{M}(B\otimes\mathcal{K})/I_{+}$ as we have done in
  proposition \ref{P:CP}. Note that \emph{the left multiplicative domain} for
  $\tilde{\phi_{+}}$ contains $z$. Hence \emph{the left multiplicative domain}
  contains $z^{n}$ for $ n \geq 0 $.
  $C(\mathbb{T})$ i.e. $\tilde{\phi_{+}}$ is homomorphism. Since every extremal partial isometry in
  $C(\mathbb{T})$ has the form $e^{ih}z^{n}$ where $h$ is
  self-adjoint element in $C(\mathbb{T})$, to finish the proof, it is enough to show that $w=\tilde{\phi}(e^{ih}z^{n})$
  is extremal partial isometry. Since $\tilde{\phi_{+}}$ is homomorphism, we can deduce $$
  \tilde{\phi_{+}}(e^{ih}z^{n})=e^{i \tilde{\phi_{+}}(h)}(\tilde{\phi_{+}}(z))^{n}$$
  It is easily shown that  $1-w^{*}w \in I_{+}$. \\
  Similarly, using $\phi^{-}$, we get $1-ww^{*} \in I_{-}$
  also. Hence $$(1-w^{*}w)
  \mathcal{M}(B\otimes\mathcal{K})(1-ww^{*})=0$$ Consequently, $w$
  is an extremal partial isometry.
\begin{corollary}
If $\phi$ is a *-homorphism from $C(\mathbb{T})$ to $B$  a
$C^{*}$-algebra, then it is extremal preserving.
\end{corollary}
\emph{proof.} It is straightforward.

By applying the definition\ref{D:cycle} to $S$ and $z$ we have
$\mathcal{E}(S,B)$ be the set of the pairs
$(\phi_{+},\phi_{-})$\,s.t.
\begin{enumerate}
\item $\phi_{+}:S \to \mathcal{M}(B\otimes\mathcal{K}) $
*-homomorphism,
      $ \phi_{-}: S \to \mathcal{M}(B\otimes\mathcal{K})$ completely positive contractive s.t.
    \begin{enumerate}
      \item $\tilde{\phi_{-}}(z)$ is extremal partial isometry for the unitary element $z \in \tilde{S}=C(\mathbb{T})$ \\
      \item $ \tilde{\phi_{-}}(z^{n})=(\tilde{\phi_{-}}(z))^{n}$
            $\tilde{\phi_{-}}((z^{*})^{n})=(\tilde{\phi_{-}}(z^{*}))^{n}$ for the unitary element $u \in \tilde{S}$
     \end{enumerate}
\item $\tilde{\phi_{+}}(z) - \tilde{\phi_{-}}(z) \in
B\otimes\mathcal{K}$ thus $(\tilde{\phi_{+}}(z))^{*}
\tilde{\phi_{-}}(z) \in 1+B\otimes\mathcal{K} $ \\
\end{enumerate}

   We call such a pair $(\phi_{+},\phi_{-})$ as \emph{an extremal cycle} and
denote the set of extremal cycles by
   $\mathcal{E}(S,B)$. Two extremal cycles
$(\phi_{+},\phi_{-}),(\psi_{+},\psi_{-})$ are \emph{homotopic}
   when there is an extremal cycle from $S$ to $\mathcal{M}([0,1]B\otimes
   \mathcal{K})$ i.e. $(\lambda_{+}^{t},\lambda_{-}^{t})\in
   \mathcal{E}(S,B), t \in [0,1] $ s.t.
   \begin{enumerate}
   \item the maps $t \to \lambda_{+}^{t}(f)+1$ and $ t \to
\lambda_{-}^{t}(f)+1$ from $[0,1]$ to $\mathcal{M}(B\otimes\mathcal{K})$
         are strictly continuous.
   \item For each $t$, $(\lambda_{+}^{t}(f)+1)^{*}, \lambda_{-}^{t}(f)+1$ are composable
   and the map $ t \to (\lambda_{+}^{t}(f)+1)^{*}\lambda_{-}^{t}(f)+1$
   from $[0,1]$ to
   $1+B\otimes\mathcal{K}$ is norm continuous.
   \item $(\lambda_{+}^{0}, \lambda_{-}^{0})=(\phi_{+},\phi_{-}),
   \quad (\lambda_{+}^{1}, \lambda_{-}^{1})=(\psi_{+},\psi_{-})$
   \end{enumerate}
   We write $(\phi_{+},\phi_{-})\sim (\psi_{+},\psi_{-})$ in this
   case.
\begin{remark}

 In fact, since
 $\mathcal{E}(S)=\mathcal{E}(\tilde{S})=\mathcal{E}(C(\mathbb{T}))=\mathcal{U}(C(\mathbb{T}))$,
 the set of homotopy classes of
$\mathcal{E}(S)$ is $K_{1}(C(\mathbb{T}))=\{[z(t)^{n}]\}\cong
\mathbb{Z}$. By proposition\ref{P:Suspension} and
corollary\ref{C:composabilityhomotopy}, for $(\phi_{+},\phi_{+})$
to be extremal cycle, it is enough to consider $f(t) \in S$ (or
$z(t)\in C(\mathbb{T}))$.
\end{remark}
  \begin{definition}
   We let $E_{\infty}(S,B)\overset{\text{def}}{=}\mathcal{E}(S,B)/\sim$ denote the
homotopy classes of extremal cycles.
  \end{definition}

Our starting point is the following proposition which tells us an
intimate relationship between KK-theory and K-theory.
\begin{lemma}\label{L:Homotopy}
  Let $w$ be a extremal partial isometry in
  $\mathcal{M}(B\otimes\mathcal{K})$. Then there is a strictly
  continuous path $w_{t}, t\in[0,1]$, of extremal partial
  isometries such that $w_{0}=1, w_{1}=w$.Furthermore, if $w$ is of the form $1+
  B$, then we can take $w(t)$ of the form $1+B\otimes \mathcal{K}$.
\end{lemma}
\emph{proof}
  Since $B\otimes\mathcal{K}$ is stable $C^{*}$-algebra, there is a path $v_{t}$ , $t \in ]0,1]$,
  of isometries in $\mathcal{M}(B\otimes\mathcal{K})$ such that
 \begin{enumerate}
 \item the map $t \mapsto v_{t}$  is strictly continuous,
 \item $ v_{1}=1$, and
 \item $ \lim_{t \to 0} v_{t}v_{t}^{*} = 0 $  in the strict topology
 \end{enumerate} See \cite{JK}.
  Set $w(t)=v_{t}wv_{t}^{*}+1-v_{t}v_{t}^{*},\, t \in ]0,1]$, and $w_{0}=1$.
  We leave the reader to check $w_{t}$ has the desired properties.

\begin{proposition} \label{P:1}
$KK(S,B) \cong K_{1}(B)$ where $S$ is suspension and B is
trivially graded stable $C^{*}$-algebra.
\end{proposition}

\emph{proof.}
  We give a proof based on Cuntz picture of KK. Recall that
\[
KK(S,B)= \{ [\phi_{+},\phi_{-}] | \phi_{+}, \phi_{-}: S \to
\mathcal{M}(B\otimes\mathcal{K}) \, \mbox{s.t.} \,
\phi_{+}-\phi_{-} \in B \otimes \mathcal{K} \}
\] For this definition, you can refer p155-156 \cite{Bl}.\\
Observe that any *-homomorphism $\phi$ from $S$ into a unital
$C^{\sp*}$-algebra defines a unitary $\phi(f)+1$ where
$f(t)=e^{2\pi it}-1$. Conversely, any unitary $u$ defines a
homomorphism by sending $f$ to $u-1$. Two homomorphisms are
homotopic if and only if the corresponding unitaries are
homotopic. From this,if we let $U_{+}, U_{-}$ be the unitaries
which come from $\phi_{+},\phi_{-}$ respectively then we get
\[
KK(S,B) \cong \{ [U_{+},U_{-}] | U_{+}-U_{-} \in B\otimes
\mathcal{K} \; U_{+}, U_{-} \in \mathcal{M}( B\otimes\mathcal{K} )
\}
\]
Note that $U_{+}^{*}U_{-}$ is a unitary in $( B\otimes \mathcal{K}
)^{\sim}$ which is a unitization of $B\otimes \mathcal{K}$. So we
can define a map $\Delta:KK(S,B) \to K_{1}(B\otimes \mathcal{K})$
by $\Delta([U_{+},U_{-}])=[U_{+}^{*}U_{-}]_{1} $. Since
$[U_{+},U_{-}]=[V_{+},V_{-}]$ implies that there exist maps $t \to
W_{\pm}^{t}$ from $[0,1]$ to $\mathcal{M}( B\otimes\mathcal{K} )$
s.t. $t \to  W_{+}^{t}-W_{-}^{t}$ is continuous in norm $
B\otimes\mathcal{K}$ and $(W_{+}^{0},W_{-}^{0})=(U_{+},U_{-}),\,
(W_{+}^{1},W_{-}^{1})=(V_{+},V_{-}) $ , $U_{+}^{*}U_{-}$ is
homotopic to $V_{+}^{*}V_{-}$. Hence $\Delta$ is  well defined. Also, you can easily check $[U_{+},U_{-}]$ are degenerate
if and only if $U_{+}=U_{-}$.\\ 
To prove $\Delta$ is surjective, let $U$ be the unitary in $(B\otimes \mathcal{K} )^{\sim}$. As we oberved, there is a unique map $\phi:S \mapsto B\otimes \mathcal{K}$ s.t. $\phi(z-1)=U-1 $. Then $\Delta([0,\phi])= U$. We remark that we could associate to $[0,\phi]$ Kasparov module $[(\hat{\mathbb{H}}_{B},\phi,0)]$ using compact perturbation.\\   
Now let $\Delta([U_{+},U_{-}])=\Delta([V_{+},V_{-}])$.
Then there is a homotopy between  $ U_{+}^{*}U_{-}$ and
$V_{+}^{*}V_{-}$. let $P(t)$ be the corresponding homotopy in $1 +
( B \otimes K )$. Since there are strictly continuous maps
$V(t):[0,1] \mapsto \mathcal{M}(B\otimes\mathcal{K})$ ,
$W(t):[0,1] \mapsto \mathcal{M}(B\otimes\mathcal{K})$ which
connects $U_{+}$ and $V_{+}$, $U_{-}$ and $V_{-}$ respectively by
lemma \ref{L:Homotopy}, then $(W(t)P(t)^{*},V(t)P(t)$ is a
homotopy between $(U_{+},U_{-})$ and $(V_{+},V_{-})$. Hence we
have proven $ \Delta $  is injective.

 The following lemmas are the key facts in this paper
which
  employ nice properties of $S$ and $\mathcal{M}(B\otimes\mathcal{K})
  $.
\begin{lemma}\label{L:Existece}
  Given an extremal partial isometry $v$ in a unital
  $C^{*}$-algebra $B$, there is a completely positive contractive
  map from $S$ to $B$ which sends $f$ to $v-1$.In fact, there is
  the unique completely positive unital map $\tilde{\phi}$ from $C(\mathbb{T})$ to $B$
  which sends $z$ to $v$ such that $\tilde{\phi}(z^{n})=v^{n}$ and $\tilde{\phi}(z^{-n})=(v^{*})^{n}$ for $n \geq 0$.
\end{lemma}
\emph{proof.}
  By lemma\ref{L:S}, it is enough to show that there is a contractive positive map
  $\phi: S \to B $ which send $f$ to $v-1$. Let $I_{+} , I_{-}$ be defect ideals  of $v$.
  Define $\tilde{\phi}$ by $\tilde{\phi}(
  p(e^{i\theta})+\overline{q(e^{i\theta})})=p(v)+q(v)^{*}$ where $p,q$ are polynomials in $C(\mathbb{T})$.
   In $B/I_{+}$, $\pi^{+}(v)=\overline{v}$ is an isometry.
 Therefore $\tilde{\phi}^{+}(e^{i\theta})=\overline{v}$.
 If $\tau(e^{i\theta})= \sum_{n=-N}^{N}a_{n}e^{in\theta}$ is a positive function in
 $C(\mathbb{T})$,
  then there is a function $f(z)= \sum_{n=0}^{N}b_{n}z^{n}$ such
  that $\tau(e^{i\theta})=|f(e^{i\theta})|^{2}$. It is easy to check
  $\tilde{\phi}^{+}(f(e^{i\theta})\overline{f(e^{i\theta})})={f(\overline{v})}f(\overline{v})^{*}$
  since $\overline{v}$ is isometry. Hence, $\tilde{\phi}^{+}$ is positive. By Russo-Dye theorem,
  $\|\tilde{\phi^{+}}\|=\|\tilde{\phi^{+}}(1)\|=1$. By
  this, $\tilde{\phi}^{+}$is also contractive. Similarly, we can show that $\tilde{\phi}^{-}$ is positive and
  contractive. Then, by the proposition\ref{P:CP}, $\tilde{\phi}$ is
  completely positive and contractive. Let $\phi$ be the
  restriction of $\tilde{\phi}$ to $S$. Then $\phi$ is also completely positive and  \[\sup_{a \in
  S}\frac{\|\phi(a)\|}{\|a\|}=\sup_{a \in
  S}\frac{\|\tilde{\phi}(a)\|}{\|a\|} \leq \sup_{a \in
  C(\mathbb{T})}\frac{\|\tilde{\phi(a)}\|}{\|a\|}\leq 1 \] implies
  it is also contractive. From the definition of $\tilde{\phi}$,
  it is unique.

\begin{lemma}\label{L:KuiMin}
The unitary group of $\mathcal{M}(B\otimes\mathcal{K})$ is
path-connected in norm topology.
\end{lemma}
\emph{proof.} This is well-known Kuiper-Mingo's theorem. The proof
of this theorem can be found in \cite{CH} or \cite{Mi}.

Now we are ready to prove the main result of this paper. By our
observation, if we let $u$, $v$ be an unitary element in
$\mathcal{M}(B\otimes\mathcal{K})$, an extremal partial isometry
in $\mathcal{M}(B\otimes\mathcal{K})$ from $\tilde{\phi_{+}}$,
$\tilde{\phi_{-}}$ respectively and denote by $[u,v]$ its homotopy
class, then we have natural map which send $[\phi_{+},\phi_{-}]$
to $[u,v]$. In fact, we have
$$ E_{\infty}(S,B)\cong \{[u,v] \mid u^{*}v \in \mathcal{E}(\mathcal{M}(B\otimes\mathcal{K}))
\, \& \,  u^{*}v \in 1+B\otimes\mathcal{K} \}$$\\
Clearly, this map is well-defined and it is surjective by the
lemma\ref{L:Existece}.
  Suppose $\tilde{\phi_{+}^{i}}(z)=u_{i}$ , $\tilde{\phi_{-}^{i}}(z)=v_{i}$ for $i=0,1$  If
$[u_{0},v_{0}]=[u_{1},v_{1}]$, there is a homotopy between
$(u_{0},v_{0})$ and $(u_{1},v_{1})$. i.e. there are strictly
continous maps $\lambda_{\pm}: [0,1]\mapsto
U(\mathcal{M}(B\otimes\mathcal{K}))$ s.t. $\lambda_{+}-\lambda_{-}
$ is norm-continuous map \&
$(\lambda_{+}(i),\lambda_{-}(i))=(u_{i},v_{i}) $ for $i=0,1$. By
lemma \ref{L:Existece}, there are corresponding maps
$\lambda_{\pm}^{t}:S \mapsto \mathcal{M}(B\otimes\mathcal{K})$ for
each $t\in [0,1]$. We must show the strict continuity of each map.
Using the fact $\{\sum_{n=-k}^{n=k}a_{n}z^{n}\ \mid a_{n} \in
\mathbb{C}\}$ is dense in $C(\mathbb{T})$, it's enough to show $t
\to \lambda_{\pm}^{t}(\sum_{n=-k}^{n=k}a_{n}z^{n})T$ where $T \in
\mathcal{M}(B\otimes\mathcal{K})$ is norm continuous with respect
to $t$. Given $\epsilon
> 0$, $ \|\lambda_{\pm}^{t}(\sum_{n=-k}^{n=k}a_{n}z^{n})T -
\lambda_{\pm}^{s}(\sum_{n=-k}^{n=k}a_{n}z^{n})T \| \leq \|
\sum_{n=-k}^{n=k}a_{n}({\lambda_{\pm}^{t}(z)}^{n}-
{\lambda_{\pm}^{s}(z)}^{n})T \| $. Now let $\delta$ be such that
if $|t-s|< \delta$ then  $
\|(\lambda_{\pm}(t))^{n}-(\lambda_{\pm}(t))^{n}T\|<
\frac{\epsilon}{2k+1sup_{n}{|a_{n}|}}$ for each $n= -k, -k+1,
\dots , k-1, k $. Therefore if $|t-s|< \delta $  we have $ \|
\sum_{n=-k}^{n=k}a_{n}({\lambda_{\pm}^{t}(z)}^{n}-
{\lambda_{\pm}^{s}(z)}^{n})T \|<
\sum_{n=-k}^{n=k}|a_{n}|\|(\lambda_{\pm}(t))^{n}-(\lambda_{\pm}(t))^{n}T
 \| <\epsilon $. \\
 Similarly, $t \to T \lambda_{\pm}^{t}(\sum_{n=-k}^{n=k}a_{n}z^{n})
$ is shown to be norm-continuous with respect to $t$ and   $t \to
\lambda_{+}^{t}(f) - \lambda_{-}^{t}(f)$ is shown to be
norm-continuous with respect to $t$ for each
$f=\sum_{n=-k}^{n=k}a_{n}z^{n}$.\\
 Finally we should check that
$(\lambda_{+}^{i},\lambda_{-}^{i})=( \phi_{+}^{i},\phi_{-}^{i}) $
for $i=0,1$. But both maps $\tilde{\lambda_{\pm}^{i}}$ and
$\tilde{\phi_{\pm}^{i}}$ send $z$ to same extremal partial
isometries $u_{i}$ and $v_{i}$. In addition, both maps send
$z^{n}$ to $u_{i}^{n}$ , $v_{i}^{n}$ and send $(z^{*})^{n}$ to
$(u_{i}^{*})^{n}$,$(v_{i}^{*})^{n}$. From the uniqueness of
lemma \ref{L:Existece}, the conclusion follows. \\

 Since $ K_{1}(B)$ is embedded in
$[\mathcal{E}_{\infty}(B)]$, we can think of a set containing
$KK(S,B)$ which extends the map $\Delta: KK(S,B) \to K_{1}(B)$. As
indicated above (proposition\ref{P:Subset}), we have the following
theorem.

\begin{theorem}
  There is  a bijection $\Delta_{e}:[\mathcal{E}_{\infty}(S,B)]\to
  [\mathcal{E}_{\infty}(B)]$ such that the following diagram is
  commutative.

  \[
    \begin{CD}
    [\mathcal{E}_{\infty}(S,B)] @>\Delta_{e}>>   [ \mathcal{E}_{\infty}(B)
    ]\\
    @AAA                    @AAA \\
    KK(S,B)@>\Delta>>      K_{1}(B)
    \end{CD}
  \]
\end{theorem}
\emph{proof.}
We define the map $\Delta_{e}$ by $\Delta_{e}([u,v])=u^{*}v$ as we have defined $\Delta$.\\
It's not hard to check well-definess of the map.(It is almost same
as the well-definess of the map $\Delta$.)\\
 {\bf Surjectivity}: Let $v$ be the extremal partial isometry in
$\widetilde{B\otimes\mathcal{K}}$. Since the quotient map from
$\widetilde{B\otimes\mathcal{K}}$ onto
$\frac{\widetilde{B\otimes\mathcal{K}}}{B\otimes\mathcal{K}}=
\mathbb{C}$ is extremal preserving map, the scalar part of v is
also the extremal partial isometry in $\mathbb{C}$. If v is
written as $\lambda + T $ where $T \in B\otimes\mathcal{K}$ and
$\lambda \in \mathbb{C}$, $(1-\left|
\lambda\right|^{2})\mathbb{C}(1-\left|\lambda\right|^{2})=0 $.
Hence $1-\left|\lambda\right|^{2}=0 $. In other words, $\lambda
\in \mathbb{T}$. Hence we may assume $v$ is of the form $1+
B\otimes \mathcal{K}$ if necessary to multiply $\bar{\lambda}$.
Then there is a map $\phi:S \to \mathcal{M}(B\otimes\mathcal{K})$
such that $\phi(f)+1=v$ by the
lemma \ref{L:Existece}. Then $\Delta_{e}([0,\phi])=[v]$. \\
{\bf Injectivity}:  Let $\Delta_{e}([u_{0},v_{0}])=
\Delta_{e}([u_{1},v_{1}])$. We may assume
$(u_{0})^{*}v_{0}=(u_{1})^{*}v_{1}$. Let $\lambda_{+}:[0,1]
\mapsto U(\mathcal{M}(B\otimes\mathcal{K}))$ be a norm-continuous
map between $u_{0}$ and $u_{1}$ by lemma \ref{L:KuiMin}. Set
$\lambda_{-}(t)= \lambda_{+}(t)((u_{0})^{*}v_{0})$.\\
Then $$ \lambda_{-}(0)=\lambda_{+}(0)(u_{0})^{*}v_{0}=v_{0}$$
     $$ \lambda_{-}(1)=\lambda_{+}(1)(u_{1})^{*}v_{1}=v_{1}$$
Also, $\lambda_{+}(t) -
\lambda_{-}(t)=\lambda_{+}(t)(1-(u_{0})^{*}v_{0}) \in
B\otimes\mathcal{K}$ for all $t\in [0,1]$. \\
Therefore $(\lambda_{+},\lambda_{-})$ is a homotopy between
$(u_{0},v_{0})$ and $(u_{1},v_{1})$ as we wanted.

\begin{theorem}
There is a bijective map $\Delta_{k}$  from $ KK_{e}(S,B) $ onto
$K_{e}(B)$ such that the following diagram is commutative.
\[
    \begin{CD}
    KK_{e}(S,B)                @>\Delta_{k}>>    K_{e}(B) \\
    @A\kappa_{e}AA                             @AA\kappa A \\
    [\mathcal{E}_{\infty}(S,B)] @>\Delta_{e}>>  [E_{\infty}(B)] \\
    \end{CD}
  \]
\end{theorem}

\emph{proof.} For each extremal cycle $(\phi_{+},\phi_{-})$, we
shall denote the element of $KK(S,B)$ by $[(\phi_{+},\phi_{-})]$ (
or equivalently, $[u,v]$ ) To avoid confusion we shall denote by
$[(\phi_{+},\phi_{-})]_{\infty}$ its homotopy class in
$\mathcal{E}_{\infty}(S,B)$. Similarly, for each $w$ in
$\mathcal{E}(\widetilde{B\otimes\mathcal{K}})$, we shall denote by
$[w]$ its equivalent class in $K_{e}(\widetilde{B\otimes
\mathcal{K}})$ and by $[w]_{\infty}$ its homotopy class in
$E_{\infty}(\widetilde{B\otimes \mathcal{K}})$. Now define
$\Delta_{k}([\phi_{+},\phi_{-}])$ by $[u^{*}v]$ where
$u=\tilde{\phi_{+}}(z)$ and $v=\tilde{\phi_{-}}(z)$.\\
Let $[(\phi_{+},\phi_{-})]= [(\psi_{+},\psi_{-})]$. Then there is
$(\tau_{+},\tau_{-})$ such that $\tilde{\tau_{-}}(z)=v $ has
smaller defects ideals than $\tilde{\phi_{-}}(z)=v_{0}$ and
$\tilde{\psi_{-}}(z)=v_{1}$ s.t.
 \[\left(  \Theta_{B}\circ \left[
                                \begin{matrix}
                                \phi_{+} & 0     \\
                                        0        & \tau_{+}
                               \end{matrix}
                           \right], \Theta_{B}\circ \left[
                                \begin{matrix}
                                \phi_{-} & 0     \\
                                        0        & \tau_{-}
                               \end{matrix}
                           \right] \right) \sim \left( \Theta_{B}\circ\left[
                                \begin{matrix}
                                \psi_{+} & 0     \\
                                        0        & \tau_{+}
                               \end{matrix}
                           \right], \Theta_{B}\circ \left[
                                \begin{matrix}
                                \psi_{-} & 0     \\
                                        0        & \tau_{-}
                               \end{matrix}
                           \right] \right) \]
  Since $\Theta_{B}$ is isomorphism, this implies
 \[  \left( \left[
                                \begin{matrix}
                                \tilde{\phi_{+}}(z) & 0     \\
                                        0        & \tilde{\tau_{+}}(z)
                               \end{matrix}
                           \right]\right)^{*} \left[
                                \begin{matrix}
                                \tilde{\phi_{-}}(z) & 0     \\
                                        0        &
                                        \tilde{\tau_{-}}(z)
                               \end{matrix}
                           \right] \sim \left(\left[
                                \begin{matrix}
                                \tilde{\psi_{+}}(z) & 0     \\
                                        0        &
                                        \tilde{\tau_{+}}(z)
                               \end{matrix}
                           \right]\right)^{*}  \left[
                                \begin{matrix}
                               \tilde{ \psi_{-}}(z) & 0     \\
                                        0        &
                                        \tilde{\tau_{-}}(z)
                               \end{matrix}
                           \right]
\]Therefore
$[u_{0}^{*}v_{0}]_{\infty}+[u^{*}v]_{\infty}=
[u_{1}^{*}v_{1}]_{\infty}+[u^{*}v]_{\infty}$ i.e.
$[u_{0}^{*}v_{0}]_{\infty} \approx [u_{1}^{*}v_{1}]_{\infty}$\\
So far we have shown $\Delta_{k}$ is well-defined.\\
From the definition of $\Delta_{k}$, the commutativity of diagram
follows easily so that the map is surjective.\\
It remains only to show the map is injective. For this let
$[u_{0}^{*}v_{0}]=[u_{1}^{*}v_{1}]$ in $K_{e}(\widetilde{B\otimes
\mathcal{K}})$. Note that defect ideals of $u_{0}^{*}v_{0}$ and
$u_{1}^{*}v_{1}$ are same to defect ideas of $v_{0}$ and $v_{1}$.
Therefore there is $v$ in $\mathcal{E}(\widetilde{B\otimes
\mathcal{K}})$ such that $v$ has smaller defect ideals that
$v_{0}$ and $v_{1}$ s.t. $u_{0}^{*}v_{0}v \sim u_{1}^{*}v_{1}v $.
In other words,
$[(u_{0},v_{0})]_{\infty}+[(1,v)]_{\infty}=[(u_{1},v_{1})]_{\infty}+[(1,v))]_{\infty}$.
This implies the map is injective  by the routine argument. \\


\begin{thebibliography}{99}

\bibitem[BrPed]{BrPe}Lawrence G. Brown, Gert K. Pedersen  \emph{Extremal
K-Theory and Index for $C^{*}$-algebras}
K-Theory 20:201-241, 2000

\bibitem[BrPed2]{BrPe2}Lawrence G. Brown, Gert K. Pedersen \emph{On the geometry of Unit ball of a $C^{*}-algebra$}
J. Reine agnew. Math. 469 (1995), 113-147

\bibitem[Bl]{Bl} Bruce Blackadar \emph{K-theory for Operator
Algebras} MSRI Publicatons Vol.5 Second Edition Cambridge
University Press 1998

\bibitem[CH]{CH} Joachim Cuntz, Nigel Higson \emph{Kuiper's theorem for Hilbert
modules}Contemp. Math. 62 (1987), 429-434

\bibitem[JenThom]{JK}  Kjeld K. Jensen, K.Thomsen \emph{Elements of KK-theory} Birkh\"{a}user, Boston, 1991 MR94b:19008

\bibitem[Kas1]{Kas80}  G.G. Kasparov \emph{Hilbert $C^{*}$-modules:
Theorem of Stinespring and Voiculescu} J.Operator Theory 4 (1980)
133-150

\bibitem[Kas2]{Kas81} G.G. Kasparov \emph{The Operator K
functor and extensions of $C^{\sp*}$algebras} Mathe. USSR-Izv
16(1981) 513-572[English Tanslation]

\bibitem[Lan]{Lan} E.C. Lance \emph{Hilbert $C^{*}$-modules A tookit for operator algebraists}
 London Mathematical Society Lecture Note Series 210 Cambride
Univerity Press

\bibitem[Mi]{Mi} J.A.Mingo \emph{K-theory and multipliers of stable
$C^{*}$-algebras} Trans.Amer.Math.Soc. 299, 1(1987)

\bibitem[P]{Pl} Vern Paulsen \emph{Completely Bounded Maps and Operator Algebras}
Cambridge Studies in Advanced Mathematics 78 2002
\end{thebibliography}
\end{document}